\documentclass[11pt,reqno]{amsart}

\usepackage[T1]{fontenc}
\usepackage{lmodern}
\usepackage{amsmath,amssymb,amsfonts,amsthm,mathtools,mathrsfs}
\usepackage[margin=1.05in]{geometry}
\usepackage{microtype}
\usepackage{enumitem}
\usepackage[colorlinks=true,linkcolor=blue!55!black,citecolor=blue!55!black,urlcolor=blue!55!black]{hyperref}
\usepackage{xcolor}

\numberwithin{equation}{section}
\theoremstyle{plain}
\newtheorem{theorem}{Theorem}[section]

\newtheorem{lemma}[theorem]{Lemma}
\newtheorem{corollary}[theorem]{Corollary}
\theoremstyle{definition}

\theoremstyle{remark}
\newtheorem{remark}[theorem]{Remark}

\DeclareMathOperator{\GL}{GL}
\DeclareMathOperator{\Ind}{Ind}
\DeclareMathOperator{\Hom}{Hom}
\DeclareMathOperator{\diag}{diag}
\DeclareMathOperator{\sgn}{sgn}

\DeclareMathOperator{\Rea}{Re}
\DeclareMathOperator{\Ima}{Im}
\newcommand{\C}{\mathbb C}
\newcommand{\R}{\mathbb R}
\newcommand{\Z}{\mathbb Z}
\newcommand{\calJ}{\mathcal J}
\newcommand{\calP}{\mathscr P}
\newcommand{\calC}{\mathscr C}
\newcommand{\abs}[1]{\lvert #1\rvert}

\newcommand{\sph}{\mathbb S}
\newcommand{\dd}{\,d}

\title[Archimedean Flicker--Rallis periods]{Exact Archimedean Flicker--Rallis Periods on Cohomological Minimal $K$-Types}
\author{Zheng WANG}

\address{LAGA UMR 7539, Institut Galil\'ee,
Universit\'e Sorbonne Paris Nord (Paris 13)}
\email{zheng.wang@math.univ-paris13.fr}


\begin{document}
\begin{abstract}
For the symmetric pair $(\mathrm{GL}_n(\mathbb C),\mathrm{GL}_n(\mathbb R))$, we compute the exact Flicker--Rallis period of a distinguished vector in a cohomological minimal $U(n)$-type, normalized by the Gel'fand--Tsetlin pairing of Hara--Miyazaki--Namikawa. Our main tool is a weight-independent identity between the Whittaker period and a compact $O(n)$-average. We determine the proportionality constant using the archimedean newform formula of Humphries--Jo, and evaluate the compact functional on the normalized distinguished vector. This yields an explicit formula for the Flicker--Rallis period for every cohomological weight.
\end{abstract}
\maketitle

\section{Introduction}\label{sec:intro}
Flicker--Rallis periods play a fundamental role in the
theory of distinguished automorphic representations
and Asai $L$-functions. Let $E/F$ be a quadratic
extension of number fields. The global Flicker--Rallis
period is obtained by integrating cuspidal automorphic
forms on $\GL_n(\mathbb A_E)$ along
$\GL_n(\mathbb A_F)$, modulo the appropriate rational
and central subgroups, possibly with a character twist.
In his pioneering work \cite{Flicker1988,Flicker1991},
Flicker established fundamental connections between
distinction and the analytic behavior of the twisted
tensor, or Asai, $L$-function, particularly its pole
at $s=1$. Further developments concerning the poles
of these $L$-functions were obtained by Flicker and
Zinoviev \cite{FlickerZinoviev1995}. The relationship
between distinguished representations and base change
from unitary groups was also investigated through
relative summation formulas in \cite{Flicker1992}.

An important aspect of this theory is the realization of global periods in terms of local Whittaker functionals. Through the unfolding of the relevant global integral representations, Flicker--Rallis periods give rise to local integrals over mirabolic quotients. In particular, Zhang \cite[Proposition~3.2]{Zhang2014} established an explicit factorization of the global Flicker--Rallis period in terms of local Whittaker integrals. This makes the choice of local test vectors essential: for arithmetic applications, one needs not only nonvanishing but also exact evaluations of these functionals on suitably normalized vectors. Such calculations arise, in particular, in the construction of $p$-adic $L$-functions; see \cite{Wang2026}.

In this paper, we focus on the nonsplit archimedean
case of the Flicker--Rallis period, corresponding to
the pair
\[
  (G_n,H_n)=(\GL_n(\C),\GL_n(\R)).
\]
Our aim is to compute the exact Flicker--Rallis period
on a certain distinguished vector in a cohomological minimal
$U(n)$-type, with all normalization choices made
explicit.

On the analytic side, Humphries developed a theory of
canonical archimedean newforms in \cite{Humphries}.
Humphries and Jo subsequently evaluated the
Flicker--Rallis periods of these newforms in terms
of explicit local $L$-factors
\cite[Theorem~6.9]{HumphriesJo}.
However, these newforms do not, in general, belong
to the cohomological minimal $U(n)$-type.
Consequently, their period formula does not directly
determine the value on a cohomologically normalized
test vector.

On the algebraic side, Hara--Miyazaki--Namikawa
\cite{HMN} constructed explicit determinantal and
Gel'fand--Tsetlin models of highest-weight
representations, together with normalized invariant
pairings and Whittaker embeddings. These constructions
allow them to normalize Eichler--Shimura cohomology
classes explicitly and to establish uniform integrality
results for critical Rankin--Selberg $L$-values.
Their approach also provides a natural framework for
the precise normalization of Eichler--Shimura morphisms,
an issue of particular importance in arithmetic
applications.
Indeed, Hara--Miyazaki--Namikawa explicitly motivate
their uniform integrality results by the study of
Kummer-type congruences and the construction of
$p$-adic $L$-functions. Their normalization methods
have also been used in subsequent explicit calculations
of archimedean modular symbols, notably by
Jin--Tian \cite{JinTian2026}. 

The purpose of this paper is to connect the
analytic period formula of Humphries--Jo with
the algebraic normalization of
Hara--Miyazaki--Namikawa. We adopt the latter
normalization for the distinguished vector in
the cohomological minimal $U(n)$-type. 
Pattanayak--Wu--Zhang \cite{PWZ}
have established nonvanishing results on
minimal $K$-types, and more recently, Chen--Cheng \cite{CC26}
used theta lifting to obtain explicit Whittaker
functions on minimal $K$-types of highest weight
$(\kappa,0,\ldots,0)$ and to evaluate the associated
Asai local zeta integrals. Our result determines the exact Flicker--Rallis
period of the Hara--Miyazaki--Namikawa-normalized
distinguished vector in the minimal $K$-type
for every cohomological weight.
The resulting formula is compatible with the
cohomological normalization arising from
the Eichler--Shimura morphism and supplies
the explicit archimedean factor needed in
the related construction of a $p$-adic
twisted adjoint $L$-function in
\cite{Wang2026}.

The key observation is that, for the even-weight principal series,
the Flicker-Rallis period and an elementary compact average belong to the
same one-dimensional relative Hom space.  
Write $\chi(z)=z/\abs z$. For
$\kappa\in(2\Z)^n$, let
\[
 I_\kappa=\Ind_{B_n(\C)}^{\GL_n(\C)}
       (\chi^{\kappa_1}\otimes\cdots\otimes\chi^{\kappa_n})
\]
be normalized induction. Define the Flicker--Rallis
period and the compact average by
\[
\begin{aligned}
\calP_\kappa(f)
  &=
  \int_{N_{n-1}(\R)\backslash\GL_{n-1}(\R)}
  \calJ_{\psi}f
  \bigl(\diag(h,1)\bigr)\,\dd\bar h,\\
\calC_\kappa(f)
  &=
  \int_{O(n)}f(k)\,\dd k.
\end{aligned}
\]
Here $\calJ_{\psi}$ denotes the
Jacquet integral with respect to a nondegenerate character $\psi$ trivial on \(N_n(\mathbb R)\), and $\dd\bar h$ is the quotient measure
induced by the Iwasawa normalization of Humphries--Jo.
The precise conventions are given in
Section~\ref{sec:notations}.
The first main theorem is the following identity.

\begin{theorem}
\label{thm:intro-comparison}
For $n\ge2$ and $\kappa\in(2\Z)^n$,
\[
  \calP_\kappa(f)=\frac{\pi^{n/2}}{\Gamma(n/2)}\,\calC_\kappa(f)
  \qquad (f\in I_\kappa).
\]
In particular, the coefficient is independent of $\kappa$.
\end{theorem}

Let $\mu=(\mu_1\ge\cdots\ge\mu_n)\in\Z^n$ and let
\[
  \lambda=\lambda_{n,\mu}=2\mu^\vee+2\rho_n
   =(n-1-2\mu_n,n-3-2\mu_{n-1},\ldots,1-n-2\mu_1).
\]
The corresponding cohomological minimal $U(n)$-type is $V_\lambda$ defined in 
\cite[(2.13)--(2.14), \S2.4.2]{HMN}.
Its restriction to $O(n)$ contains a unique copy of $\det^{n-1}$.
Let $S$ denote the vector in that line normalized to pair to $1$
with HMN's unit highest-weight vector of the dual-weight representation,
and let $f_{B_n,\lambda,0}$ be the normalized embedding into $I_\lambda$
\cite[(2.48)]{HMN}. Their precise definitions are given in
Section~\ref{sec:HMN}. Define the Flicker--Rallis period by
\begin{equation}\label{eq:intro-J}
 J_n(\mu)=
 \int_{N_{n-1}(\R)\backslash\GL_{n-1}(\R)}
 \calJ_{\psi}\bigl(f_{B_n,\lambda,0}(S)\bigr)
          \bigl(\diag(h,1)\bigr)
 \sgn(\det h)^{n-1}\dd\bar h.
\end{equation}
The second theorem is the explicit value of $J_n(\mu)$.

\begin{theorem}
\label{thm:intro-main}
Let $n\ge2$ and let $\mu=(\mu_1\ge\cdots\ge\mu_n)\in\Z^n$.
Then
\[
  J_n(\mu)
    =\frac{\pi^{n/2}}{\Gamma(n/2)}.
\]
In particular, $J_n(\mu)$ is independent of $\mu$.
\end{theorem}

\subsection*{Organization of the paper.}
Section~\ref{sec:notations} fixes the notation and normalization conventions.
In Section~\ref{sec:comparison}, we establish a weight-independent identity
between the Flicker--Rallis period and a compact
$O(n)$-average, using the archimedean newform formula
of Humphries--Jo.
In Section~\ref{sec:HMN}, we evaluate this identity on the
Hara--Miyazaki--Namikawa-normalized distinguished vector
in the cohomological minimal $K$-type, and extend
the result to essentially tempered cohomological
representations by a determinant twist.

\section*{Acknowledgements}
The author is grateful to his advisor Jacques Tilouine for his guidance and for bringing \cite{HMN} to his attention. He also thanks Shih-Yu Chen and Daniel Disegni for helpful discussions and for encouraging him to investigate this problem. The calculation presented here provides the final missing archimedean ingredient needed to complete \cite{Wang2026}.

The author is a doctoral fellow of the MathPhDInFrance program, coordinated by the
Fondation Sciences Math{\'e}matiques de Paris (FSMP), enrolled at the
{\'E}cole Doctorale Galil{\'e}e and hosted by LAGA.

\medskip
\noindent\textit{Tool and computational resource disclosure.} AI was used to proofread a draft of the manuscript.

\medskip
\noindent\textit{Funding.}
This project has received funding from the European Union's Horizon Europe research
and innovation programme under the Marie Sk\l odowska-Curie grant agreement
No.~101126554. This work received support from the French government under the
France 2030 investment plan, as part of the Initiative d'Excellence
d'Aix-Marseille Universit\'e -- AMIDEX AMX-23-COFU-01.

\begin{center}
\includegraphics[width=2.5cm]{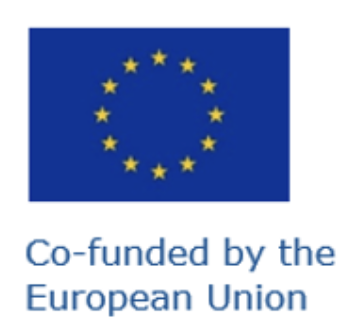}
\end{center}

\section{Basic notations}
\label{sec:notations}

Throughout, $n\ge2$, $G_n=\GL_n(\C)$, $H_n=\GL_n(\R)$, and
$K_n=U(n)$.  Write $B_n$ and $N_n$ for the upper Borel and its
unipotent radical.  Our right-translation convention for normalized
induction is
\[
 I_\kappa=\Ind_{B_n(\C)}^{G_n}
 (\chi^{\kappa_1}\otimes\cdots\otimes\chi^{\kappa_n}) =\{f\colon G_n \to \C\mid f(bg)=\delta_{B_n(\C)}^{1/2}(b)
       \prod_{j=1}^n\chi(b_{jj})^{\kappa_j}\,f(g)\}.
\]
Here $\abs z_\C=z\bar z$.  Since the inducing characters are unitary,
$I_\kappa$ is a unitary principal series; in the cases under consideration
it is irreducible and generic.  In the cohomological case its inducing
weight $\lambda_{n,\mu}$ has strictly decreasing entries, all congruent
to $n-1$ modulo $2$.

We use two nondegenerate characters of $N_n(\C)$:
\[
 \begin{split}
 \psi_+(u)&=\exp\Bigl(4\pi i\sum_{j=1}^{n-1}
                         \Rea(u_{j,j+1})\Bigr),\\
 \psi(u)&=\exp\Bigl(-4\pi i\sum_{j=1}^{n-1}
                         \Ima(u_{j,j+1})\Bigr).
 \end{split}
\]
Thus $\psi$ is trivial on $N_n(\R)$.
We put $d_\C z=2\,d\Rea z\,d\Ima z$ and use its product over
upper-unipotent coordinates.  With $w_n$ the antidiagonal matrix with all
nonzero entries equal to $1$, define the Jacquet map
\[
 \calJ_\psi f(g)=\int_{N_n(\C)}f(w_nug)\psi(u)^{-1}\dd u,
\]
initially in the region of absolute convergence
and then by holomorphic continuation in the
inducing parameters by
\cite[Theorem~15.4.1]{Wal92}; see also \cite[\S9.1]{Humphries}.

We use the same Haar measure convention as in \cite{HumphriesJo}. For clarity, we state here in full.  In the Iwasawa decomposition for $\GL_q(\R)$ we use the diagonal
$A_q(\R)=(\R^\times)^q$, the product multiplicative measure
$\prod_{j=1}^q da_j/\abs{a_j}$, the coordinate Lebesgue measure on
$N_q(\R)$, and probability Haar measure on $O(q)$:
\[
 \int_{\GL_q(\R)}F(g)\dd g
  =\int_{N_q(\R)}\int_{A_q(\R)}\int_{O(q)}
 F(uak)\,\delta_{B_q(\R)}(a)^{-1}
 \dd k\prod_{j=1}^q\frac{da_j}{\abs{a_j}}\dd u.
\]
The quotient measure on $N_q(\R)\backslash\GL_q(\R)$ is denoted
$d\bar h$.

\section{Comparison of two functionals}
\label{sec:comparison}

We first consider principal series with even inducing
exponents, for which the period is untwisted. 
The general case follows by
a determinant twist and will be treated in the end of this section.

\medskip
Let
$
 \Gamma_\R(s)=\pi^{-s/2}\Gamma(s/2),
 \Gamma_\C(s)=2(2\pi)^{-s}\Gamma(s).
$
For $\kappa\in(2\Z)^n$, set
\[
 \begin{split}
 \calP_\kappa(f)
 &=\int_{N_{n-1}(\R)\backslash\GL_{n-1}(\R)}
   \calJ_{\psi}f(\diag(h,1))\dd\bar h,\\
 \calC_\kappa(f)&=\int_{O(n)}f(k)\dd k.
 \end{split}
\]
By \cite[Lemma~6.8]{HumphriesJo}, the Flicker--Rallis
integral converges absolutely. The resulting
functional is continuous and 
mirabolic-invariant.
On the other hand, Lemma~\ref{lem:compact-invariance} and the nonvanishing
of $\calC_\kappa$, verified below in \eqref{eq:C-newform},
show that $I_\kappa$ is $H_n$-distinguished.
Thus $\calP_\kappa$ is $H_n$-invariant by
Kemarsky's theorem \cite[Theorem 1.1]{Kem15}; see also \cite[Remark~6.11]{HumphriesJo}. We now prove that $\calC_\kappa$ enjoys the same $H_n$-invariance property.
\begin{lemma}
\label{lem:compact-invariance}
For every $\kappa\in(2\Z)^n$ and every $h\in H_n$,
\[
 \calC_\kappa(\pi(h)f)=\calC_\kappa(f).
\]
Thus $\calC_\kappa\in\Hom_{H_n}(I_\kappa,\mathbf 1)$.
\end{lemma}
\begin{proof}
Let $B_\R=B_n(\R)$ and $K=O(n)$.  For $b\in B_\R$, we have
\[
 f(bh)=\delta_{B_n(\C)}^{1/2}(b) f(h)
       =\delta_{B_\R}(b) f(h).
\]
Indeed $\delta_{B_n(\C)}^{1/2}|_{B_\R}=\delta_{B_\R}$ and
$\chi^{\kappa_j}|_{\R^\times}=1$.  For $kh=b_h(k)k_h(k)$ in Iwasawa decomposition, we have
\[
 \int_K f(k_h(k))\delta_{B_\R}(b_h(k))\dd k
     =\int_K f(k)\dd k,
\]
by the Iwasawa integration
formula and the transformation of left Haar
measure on $B_{\mathbb R}$ under right translation; see for example \cite[Lemma~2.2]{Deitmar2006}. Then we conclude by
$f(kh)=\delta_{B_\R}(b_h(k))f(k_h(k))$.
\end{proof}

By the multiplicity-one
theorem for the Gel'fand pair
$(\GL_n(\C),\GL_n(\R))$~\cite{AizenbudGourevitch}, we know $\Hom_{H_n}(I_\kappa,\mathbf 1)$ is at most one-dimensional, thus $\mathcal{P}_\kappa$ is proportionate to $\mathcal{C}_\kappa$. The rest of this section is devoted to the calculation of this ratio. The most important ingredient is the archimedean newform theory introduced in \cite{Humphries} which we now recall.
\medskip

For $\kappa=(\kappa_1,\ldots,\kappa_n)\in(2\Z)^n$,
put
\[
 a_j=\frac{|\kappa_j|}{2},\qquad
 \kappa_j^+=\max(\kappa_j,0),\qquad
 \kappa_j^-=\max(-\kappa_j,0).
\]
Let $f_\kappa^\circ\in I_\kappa$ be the canonically
normalized induced newform of Humphries
\cite[Corollary~8.17]{Humphries}.
Explicitly, its restriction to $U(n)$ is 
\begin{equation} 
 f_\kappa^\circ(k)
 =
 i^{\sum_{j=1}^n(j-1)|\kappa_j|}
 \prod_{1\le i<j\le n}
       \Gamma_\C(1+a_i+a_j)
 \prod_{j=1}^n
       k_{j,n}^{\kappa_j^+}
       \overline{k_{j,n}}^{\,\kappa_j^-}.\label{eq:newform}
\end{equation}

Following
\cite[Proposition~4.2]{HumphriesJo},
let $\pi_{\mathrm{ur}}$ be the spherical
unramified companion of $I_\kappa$.
Thus, choosing a permutation $\sigma$ such that
$a_{\sigma(1)}\ge\cdots\ge a_{\sigma(n)}$, we have
\[
 \pi_{\mathrm{ur}}
 =
 \Ind_{B_n(\C)}^{\GL_n(\C)}
 \bigl(
 |\cdot|_\C^{a_{\sigma(1)}}\otimes\cdots
 \otimes|\cdot|_\C^{a_{\sigma(n)}}
 \bigr).
\]
Let $f_{\mathrm{ur}}^\circ$ be its canonically
normalized spherical section, obtained by applying
\cite[Definition~9.2]{Humphries}
to $\pi_{\mathrm{ur}}$. In particular,
\[
 f_{\mathrm{ur}}^\circ(k)
 =
 \prod_{1\le i<j\le n}
 \Gamma_\C
 \bigl(1+a_{\sigma(i)}-a_{\sigma(j)}\bigr),
 \qquad k\in U(n).
 \]
By \cite[Lemma~4.6]{HumphriesJo}, their Whittaker realizations with
respect to $\psi_+$ satisfy
\[
 \calJ_{\psi_+}f_\kappa^\circ(\diag(g,1))
 =
 \calJ_{\psi_+}f_{\mathrm{ur}}^\circ(\diag(g,1)),
 \qquad g\in\GL_{n-1}(\C).
\]
 For our later use, we also need a $\psi$-version of this comparison. We only compare them for diagonal elements.

\begin{lemma}
\label{lem:newform-transport}
For $\kappa \in (2\Z)^n$, put $\delta_\kappa=i^{\sum_{j=1}^n(j-1)\kappa_j}$. For
$a=\diag(a_1,\ldots,a_{n-1},1)$ with $a_j\in \C^\times$, 
\[
 \calJ_{\psi}f_\kappa^\circ(a)
 =\delta_\kappa^{-1}\calJ_{\psi_+}f_\kappa^\circ(a)
 =\delta_\kappa^{-1}\calJ_{\psi}
                                      f_{\mathrm{ur}}^\circ(a).
\]
\end{lemma}
\begin{proof}
Put $d_n=\diag(i^{n-1},i^{n-2},\ldots,i,1)$ and $u=d_nu'd_n^{-1}$, we have $\psi_+(u)=\psi(u')$. By the description above, the newform $f_\kappa^\circ$ is right invariant under
$\diag(U(n-1),1)$. Since $d_n$ belongs
to this subgroup and commutes with $a$, we know $
\mathcal J_{\psi_+}f_\kappa^\circ(d_na)
=\mathcal J_{\psi_+}f_\kappa^\circ(ad_n)
=\mathcal J_{\psi_+}f_\kappa^\circ(a).
$ On the other hand, for $g\in G_n$, we have
\[
 \begin{aligned}
 f_\kappa^\circ(w_nu d_ng)
  =f_\kappa^\circ(w_nd_nu'g)
  =\prod_{j=1}^n\chi(i^{j-1})^{\kappa_j}\,f_\kappa^\circ(w_nu'g)
    =\delta_\kappa f_\kappa^\circ(w_nu'g).
 \end{aligned}
\]
Thus $
\mathcal J_{\psi_+}f_\kappa^\circ(d_na)
=
\delta_\kappa\,
\mathcal J_{\psi}f_\kappa^\circ(a)
$. This proves the first equality. A similar argument shows $\calJ_{\psi_+}f_\mathrm{ur}^\circ(a)=\calJ_{\psi}f_{\mathrm{ur}}^\circ(a)$. Next,
\cite[Lemma~4.6]{HumphriesJo} gives $\calJ_{\psi_+}f_\kappa^\circ(a)=\calJ_{\psi_+}f_{\mathrm{ur}}^\circ(a) $. This gives the second equality.
\end{proof}

\begin{theorem}
\label{thm:comparison}
For every $\kappa\in(2\Z)^n$ and every $f\in I_\kappa$,
\[
 \calP_\kappa(f)=\frac{\pi^{n/2}}{\Gamma(n/2)}\calC_\kappa(f).
\]
In particular $\calC_\kappa$ and $\calP_\kappa$ are nonzero.
\end{theorem}
\begin{proof}
It suffices to prove $\calP_\kappa(f_\kappa^\circ)=\frac{\pi^{n/2}}{\Gamma(n/2)}\calC_\kappa(f_\kappa^\circ)$. We first calculate $\calC_\kappa(f_\kappa^\circ)$.

By \eqref{eq:newform}, the integral concerns only the last column whose entries are real. 
Put
$
 T=\sum_{j=1}^n a_j,
$
 using $O(n)/O(n-1)\cong \mathbb{S}^{n-1}$, the integral reduces to 
 \[
 \int_{\sph^{n-1}}\prod_{j=1}^n x_j^{2a_j}\dd\sigma(x)
 =\frac{\Gamma(n/2)}{\Gamma(n/2+T)}
   \prod_{j=1}^n\frac{\Gamma(a_j+1/2)}{\Gamma(1/2)}.
\]
One way to verify this integral is to calculate
$\int_{\R^n}\prod_j x_j^{2a_j}e^{-\sum_jx_j^2}\dd x$
using both Cartesian and polar coordinates.  We then obtain \begin{equation}\label{eq:C-newform}
 \calC_\kappa(f_\kappa^\circ)
 =\epsilon_\kappa G_\kappa
  \frac{\Gamma(n/2)}{\Gamma(n/2+T)}
  \prod_{j=1}^n\frac{\Gamma(a_j+1/2)}{\sqrt\pi}\ne0.
\end{equation}
 Here $\epsilon_\kappa=i^{\sum_{j=1}^n(j-1)|\kappa_j|},
 G_\kappa=
 \prod_{1\le i<j\le n}
 \Gamma_\C(1+a_i+a_j).$
 
To calculate $\calP_\kappa(f_\kappa^\circ)$, note that the integral reduces to positive
real Iwasawa diagonal coordinates because $f_\kappa^\circ$ is
right $O(n-1)$-invariant and $\psi$ is trivial on
$N_{n-1}(\R)$. We can thus use
Lemma~\ref{lem:newform-transport}, and then apply the archimedean
Flicker--Rallis newform formula
\cite[Theorem~6.9]{HumphriesJo}.   This gives
\[
 \mathscr P_\kappa(f_\kappa^\circ)
 =
 \delta_\kappa^{-1}
 \frac{L(1,\pi_{\mathrm{ur}},\mathrm{As})}
      {L(n,\omega_{\pi_{\mathrm{ur}}}|_{\R^\times})}.
\]
By the explicit spherical Asai $L$-factor formula
\cite[Theorem~6.3]{HumphriesJo}, we have
\[
 L(s,\pi_{\mathrm{ur}},\mathrm{As})
 =
 \prod_{j=1}^n\Gamma_\R(s+2a_j)
 \prod_{1\le i<j\le n}
       \Gamma_\C(s+a_i+a_j).
\]
On the other hand, the central character of
$\pi_{\mathrm{ur}}$ restricts to
$
 \omega_{\pi_{\mathrm{ur}}}|_{\R^\times}
 =|\cdot|^{2T},
$
and therefore
\[
 L(n,\omega_{\pi_{\mathrm{ur}}}|_{\R^\times})
 =\Gamma_\R(n+2T).
\]
Substituting these expressions into the preceding
period identity, we obtain
\begin{equation}\label{eq:P-newform}
 \mathscr P_\kappa(f_\kappa^\circ)
 =
 \delta_\kappa^{-1}
 \frac{
   \displaystyle\prod_{j=1}^n
     \Gamma_\R(1+2a_j)
   \displaystyle\prod_{1\le i<j\le n}
     \Gamma_\C(1+a_i+a_j)
 }{
   \Gamma_\R(n+2T)
 }.
\end{equation}

Finally, note that
\[
 \frac{\prod_{j=1}^n\Gamma_\R(1+2a_j)}{\Gamma_\R(n+2T)}
 =\frac{\prod_{j=1}^n\Gamma(a_j+1/2)}{\Gamma(n/2+T)}.
\]
Dividing \eqref{eq:P-newform} by \eqref{eq:C-newform},
$G_\kappa$ and all weight-dependent Gamma factors cancel. As for each $j$, $\abs{\kappa_j}-\kappa_j$ is divisible by $4$, we have
$\delta_\kappa^{-1}\epsilon_\kappa^{-1}=
\epsilon_\kappa^{-2}=1$. Thus
\[
 \frac{\calP_\kappa(f_\kappa^\circ)}
      {\calC_\kappa(f_\kappa^\circ)}
   =\frac{\pi^{n/2}}{\Gamma(n/2)}.
\]
Combining this with multiplicity one proves the assertion.
\end{proof}

We now treat the general case, allowing both parities.
Let $n\geq 2$, and let
$\lambda=(\lambda_1,\ldots,\lambda_n)\in\Z^n$
satisfy
$
 \lambda_i\equiv\varepsilon\pmod 2$ for $ 1\leq i\leq n
$
and some $\varepsilon\in\{0,1\}$.
Define
\[
 \begin{split}
 \calP_\lambda^\varepsilon(f)
 &=\int_{N_{n-1}(\R)\backslash\GL_{n-1}(\R)}
   \calJ_{\psi}f(\diag(h,1))
   \sgn(\det h)^\varepsilon\dd\bar h,\\
 \calC_\lambda^\varepsilon(f)
 &=\int_{O(n)}f(k)\det(k)^\varepsilon\dd k.
 \end{split}
\]
In particular, for the cohomological weight
$\lambda=2\mu^\vee+2\rho_n$, we have
$\varepsilon\equiv n-1\pmod 2$.

\begin{corollary}\label{cor:twist}
For the above $\lambda$ and $f\in I_\lambda$,
\[
 \calP_\lambda^\varepsilon(f)
 =(-1)^{\varepsilon\binom n2}
    \frac{\pi^{n/2}}{\Gamma(n/2)}
    \calC_\lambda^\varepsilon(f).
\]
\end{corollary}
\begin{proof}
The case $\varepsilon=0$ is
Theorem~\ref{thm:comparison}.  For $\varepsilon=1$, set
$\kappa=\lambda+(1,\ldots,1)$ and
$\tilde{f}(g)=\chi(\det g)f(g)\in I_\kappa$.
Since $\det w_n=(-1)^{\binom n2}$ and $\det u=1$ on $N_n(\C)$,
\[
 \calJ_\psi(\tilde{f})(g)
   =(-1)^{\binom n2}\chi(\det g)\calJ_\psi f(g).
\]
For real $h$, $\chi(\det h)=\sgn(\det h)$, whereas
$\chi(\det k)=\det k$ for $k\in O(n)$.
Consequently
$\calP_\lambda^1(f)=(-1)^{\binom n2}\calP_\kappa(\tilde{f})$ and
$\calC_\lambda^1(f)=\calC_\kappa(\tilde{f})$.
Then the assertion follows from Theorem~\ref{thm:comparison}.
\end{proof}

\section{Hara--Miyazaki--Namikawa normalization and the explicit Flicker-Rallis period}\label{sec:HMN}

We recall the conventions for highest-weight representations,
Gel'fand--Tsetlin bases, and invariant pairings from
\cite[\S2.3.1--2.3.3]{HMN}.

Let $\lambda=(\lambda_1,\ldots,\lambda_n)$ be a dominant
integral weight, and put
$
 \lambda^\vee=(-\lambda_n,\ldots,-\lambda_1).
$
We use the determinantal realization
$(\tau_\lambda,V_\lambda)$ of \cite[(2.13)]{HMN},
where
\[
 \tau_\lambda(g)v(z)=v(zg),
 \qquad g\in\GL_n(\C).
\]

Let $\mathrm G(\lambda)$ denote the set of integral
Gel'fand--Tsetlin patterns of type $\lambda$, as in
\cite[(2.15)]{HMN}. Thus an element
\[
 M=(m_{i,j})_{1\le i\le j\le n}
 \in\mathrm G(\lambda)
\]
is a triangular array with top row
$(m_{1,n},\ldots,m_{n,n})=\lambda$ satisfying
the interlacing conditions
\[
 m_{i,j+1}\ge m_{i,j}\ge m_{i+1,j+1}.
\]
Its weight is
$\gamma^M=(\gamma_1^M,\ldots,\gamma_n^M)$, where
\[
 \gamma_j^M
 =\sum_{i=1}^j m_{i,j}
  -\sum_{i=1}^{j-1}m_{i,j-1}.
\]
Following \cite[(2.17)--(2.21)]{HMN}, let
$\{\zeta_M\}_{M\in\mathrm G(\lambda)}$ be the
orthonormal Gel'fand--Tsetlin basis of $V_\lambda$
with respect to the $U(n)$-invariant Hermitian
pairing $(\cdot,\cdot)_\lambda$. We also introduce the rational Gel'fand--Tsetlin basis
\[
 \xi_M=\sqrt{r(M)}\,\zeta_M,
\]
where
\[
 r(M)=
 \prod_{1\le i\le j<k\le n}
 \frac{
 (m_{i,k}-m_{j,k-1}-i+j)!
 (m_{i,k-1}-m_{j+1,k}-i+j)!
 }{
 (m_{i,k-1}-m_{j,k-1}-i+j)!
 (m_{i,k}-m_{j+1,k}-i+j)!
 }.
\]
In particular,
\[
 (\xi_M,\xi_N)_\lambda
 =\delta_{M,N}r(M).
\]
For an extremal weight $\gamma\in S_n\lambda$,
let $H(\gamma)$ be the unique Gel'fand--Tsetlin pattern of weight
$\gamma$, defined above \cite[Lemma~2.2]{HMN}.
By \cite[Lemma~2.2]{HMN}, the extremal Gel'fand--Tsetlin vectors satisfy
\[
 \xi_{H(\gamma)}
 =\zeta_{H(\gamma)}
 =f_{h(\gamma)}(z).
\]
In particular, $r(H(\gamma))=1$. Let $w_0$ be the longest element of the Weyl group
$S_n$. We write
\[
 \begin{aligned}
 E_\lambda
   &:=\xi_{H(\lambda)}
     =\zeta_{H(\lambda)},\\
 L_\lambda
   &:=\xi_{H(w_0\lambda)}
     =\zeta_{H(w_0\lambda)}
 \end{aligned}
\]
for the unit highest and lowest extremal Gel'fand--Tsetlin vectors
of $V_\lambda$, respectively. We also write
\[
 E_{\lambda^\vee}
 :=\xi_{H(\lambda^\vee)}
 \in V_{\lambda^\vee}.
\]
\cite[(2.13)--(2.14) and Lemma~2.2]{HMN} gives 
an explicit determinantal description of these extremal vectors
\begin{equation}\label{eq:HMN-extremal-minors}
 \begin{aligned}
 E_\lambda(z)
 &= (\det z)^{\lambda_n}
    \prod_{j=1}^{n-1}
    \bigl(\det z_{[1,j],[1,j]}\bigr)^{
      \lambda_j-\lambda_{j+1}},\\
 L_\lambda(z)
 &= (\det z)^{\lambda_n}
    \prod_{j=1}^{n-1}
    \bigl(\det z_{[1,j],[n-j+1,n]}\bigr)^{
      \lambda_j-\lambda_{j+1}}.
 \end{aligned}
\end{equation}
Here $z_{I,J}$ is the corresponding
submatrix of $z$ for subsets $I,J\subseteq\{1,\ldots,n\}$ of the same
cardinality.
\medskip
\\
For $M\in\mathrm G(\lambda)$, define
\begin{equation}\label{eq:HMN-q}
 q(M)
 :=\sum_{j=1}^{n-1}(n-j)\gamma_j^M
 =\sum_{1\le i\le j\le n-1}m_{i,j},
\end{equation}
as in \cite[(2.28)]{HMN}.
The dual Gel'fand--Tsetlin pattern $M^\vee$ is defined by
\[
 m^\vee_{i,j}=-m_{j+1-i,j}.
\]
Thus $M^\vee\in\mathrm G(\lambda^\vee)$, and
$\gamma^{M^\vee}=-\gamma^M$. Then there is  a $\GL_n(\C)$-invariant bilinear pairing
\[
 \langle\cdot,\cdot\rangle_\lambda:
 V_\lambda\otimes V_{\lambda^\vee}\longrightarrow\C
\]
by
\begin{equation}\label{eq:HMN-Gel'fand--Tsetlin-pairing}
 \langle\xi_M,\xi_N\rangle_\lambda
 =
 \begin{cases}
  (-1)^{q(M)}r(M),&N=M^\vee,\\
  0,&N\ne M^\vee,
 \end{cases}
\end{equation}
for $M\in\mathrm G(\lambda)$ and
$N\in\mathrm G(\lambda^\vee)$; see
\cite[(2.33) and Lemma~2.9]{HMN}.

Finally, let $f_{B_n,\lambda,0}$ denote the
normalized $U(n)$-equivariant embedding of
$V_\lambda$ into the corresponding principal series,
as defined in \cite[(2.48)]{HMN}. In particular, its restriction to $U(n)$ is given by
\begin{equation}\label{eq:HMN-section-normalization}
 f_{B_n,\lambda,0}(v)(k)
 =
 (\tau_\lambda(k)v,E_\lambda)_\lambda,
 \qquad k\in U(n).
\end{equation}

We can now introduce the distinguished
vector specified by the pairing $\langle \cdot,\cdot\rangle_\lambda$ and evaluate the compact functional
$\calC_\lambda^\varepsilon$ on it, which is of most interest.

\begin{lemma}
\label{lem:HMNcompact}
Suppose $\lambda$ is dominant integral and
$
 \lambda_j\equiv\varepsilon\pmod2$
 for $1\le j\le n$ and some
 $\varepsilon\in\{0,1\}.$ Let $w_n=(\delta_{i,n+1-j})_{ij}$ be the antidiagonal matrix we fixed above. Then
the $\det^\varepsilon$-isotypic line in
$V_\lambda|_{O(n)}$ is one-dimensional, and
there is a unique vector $S_{\lambda}$
in this line such that
\[
 \langle S_{\lambda},
         E_{\lambda^\vee}\rangle_\lambda=1.
\]
Moreover,
\begin{equation}\label{eq:HMN-section-value}
 f_{B_n,\lambda,0}(S_{\lambda})(k)
 =
 (-1)^{\varepsilon {\binom n2}}\det(k)^\varepsilon,
 \qquad k\in O(n).
\end{equation}
Consequently,
\begin{equation}\label{eq:HMN-average}
 \calC_\lambda^\varepsilon
 \bigl(f_{B_n,\lambda,0}
       (S_{\lambda})\bigr)
 =
 (-1)^{\varepsilon {\binom n2}}.
\end{equation}
\end{lemma}

\begin{proof}
We first establish two identities involving
the extremal vectors.

By \eqref{eq:HMN-extremal-minors}, we obtain
\[
 \tau_\lambda(w_n)E_\lambda
 =
 (-1)^{
 {\binom n2}\lambda_n+
 \sum_{j=1}^{n-1}
 \binom j2(\lambda_j-\lambda_{j+1})
 }L_\lambda.
\]
Since right
multiplication by $w_n$ reverses the selected
columns in every minor and
$
 \det(w_n)=(-1)^{{\binom n2}},
$ Note that all differences $\lambda_j-\lambda_{j+1}$
are even, while $\lambda_n\equiv\varepsilon
\pmod2$. Hence
\[
 \tau_\lambda(w_n)E_\lambda
 =(-1)^{\varepsilon {\binom n2}}L_\lambda,
\]

Next, put
$
 M^-:=H(w_0\lambda).
$
The $j$-th row of $M^-$ is
$
 (\lambda_{n-j+1},\ldots,\lambda_n).
$
Therefore, by \eqref{eq:HMN-q},
\[
 \begin{aligned}
 q(M^-)
 =\sum_{j=1}^{n-1}
   \sum_{i=n-j+1}^n\lambda_i
 =\sum_{j=1}^n(j-1)\lambda_j
 \equiv\varepsilon {\binom n2}\pmod2.
 \end{aligned}
\]
Moreover, by\cite[Lemma 2.2]{HMN}, we know $r(M^-)=1$, and
$
 (M^-)^\vee=H(\lambda^\vee).
$
Thus \eqref{eq:HMN-Gel'fand--Tsetlin-pairing} yields
\[
 \langle L_\lambda,E_{\lambda^\vee}\rangle_\lambda
 =(-1)^{q(M^-)}r(M^-)
 =(-1)^{\varepsilon {\binom n2}}.
\]

We now justify the normalization of the
distinguished vector. By the Cartan--Helgason
theorem \cite[Theorem~3.1 and
Corollary~3.1]{He07}, the
$\det^\varepsilon$-isotypic subspace of
$V_\lambda|_{O(n)}$ is one-dimensional.
To see that its pairing with
$E_{\lambda^\vee}$ is nonzero, consider the
functional on $V_{\lambda^\vee}$ defined by
\[
\mathcal{C}(v)
 =
 \int_{O(n)}
 v(k)\det(k)^\varepsilon\,dk.
\]
We have $
\mathcal C\in
\operatorname{Hom}_{O(n)}
\bigl(V_{\lambda^\vee},\det^\varepsilon\bigr)$.
Since $\lambda^\vee_j\equiv\varepsilon\pmod2$,
\eqref{eq:HMN-extremal-minors} for
$E_{\lambda^\vee}$ gives
\[
 \begin{aligned}
\mathcal{C}(E_{\lambda^\vee})
 &=
 \int_{O(n)}
 (\det k)^{\varepsilon+\lambda^\vee_n}
 \prod_{j=1}^{n-1}
 \bigl(\det k_{[1,j],[1,j]}\bigr)^{
  \lambda^\vee_j-\lambda^\vee_{j+1}}
 \,dk\\
 &=
 \int_{O(n)}
 \prod_{j=1}^{n-1}
 \bigl(\det k_{[1,j],[1,j]}\bigr)^{
  \lambda^\vee_j-\lambda^\vee_{j+1}}
 \,dk
 >0.
 \end{aligned}
\]
Indeed, the integrand is nonnegative and
strictly positive in a neighborhood of the
identity.

By the nondegenerate invariant pairing
$\langle\cdot,\cdot\rangle_\lambda$, the
functional $\mathcal{C}$ corresponds to a
nonzero vector in the
$\det^\varepsilon$-isotypic line of
$V_\lambda$. The preceding positivity
shows that this vector pairs nontrivially
with $E_{\lambda^\vee}$. Consequently,
$S_{\lambda}$ is well-defined.

We now determine its normalization with
respect to the Hermitian pairing.
By \eqref{eq:HMN-Gel'fand--Tsetlin-pairing} and the
orthogonality of the Gel'fand--Tsetlin basis, we have
\[
 \langle v,E_{\lambda^\vee}\rangle_\lambda
 =
 (-1)^{\varepsilon {\binom n2}}\,(v,L_\lambda)_\lambda,
 \qquad v\in V_\lambda.
\]
In particular, the normalization of $S_{\lambda}$
implies
\begin{equation}\label{eq:S-lowest-pairing}
 1=(-1)^{\varepsilon {\binom n2}}\,(S_{\lambda},L_\lambda)_\lambda.
\end{equation}

Since $S_{\lambda}$ belongs to the
$\det^\varepsilon$-isotypic line,
$
 \tau_\lambda(w_n)S_{\lambda}
 =\det(w_n)^\varepsilon S_{\lambda}
 =(-1)^{\varepsilon {\binom n2}}S_{\lambda}.
$
By the $U(n)$-invariance of the Hermitian
pairing $(\cdot,\cdot)_\lambda$, we obtain
\[
 \begin{aligned}
 (-1)^{\varepsilon {\binom n2}}(S_{\lambda},L_\lambda)_\lambda
 =(S_{\lambda},\tau_\lambda(w_n)E_\lambda)_\lambda
 =(\tau_\lambda(w_n)^{-1}S_{\lambda},E_\lambda)_\lambda
 =(-1)^{\varepsilon {\binom n2}}(S_{\lambda},E_\lambda)_\lambda.
 \end{aligned}
\]
Thus
$
 (S_{\lambda},L_\lambda)_\lambda
 =(S_{\lambda},E_\lambda)_\lambda.
$
Together with \eqref{eq:S-lowest-pairing},
this gives
\[
 (S_{\lambda},E_\lambda)_\lambda=(-1)^{\varepsilon {\binom n2}}.
\]

Finally, by \eqref{eq:HMN-section-normalization},
for every $k\in O(n)$ we have
\[
 f_{B_n,\lambda,0}(S_{\lambda})(k)
 =(\tau_\lambda(k)S_{\lambda},E_\lambda)_\lambda
 =\det(k)^\varepsilon
   (S_{\lambda},E_\lambda)_\lambda
 =(-1)^{\varepsilon {\binom n2}}\det(k)^\varepsilon.
\]
This proves \eqref{eq:HMN-section-value}.
Integrating with $\det(k)^\varepsilon$ gives
\[
 \calC_\lambda^\varepsilon
 \bigl(f_{B_n,\lambda,0}(S_{\lambda})\bigr)
 =(-1)^{\varepsilon {\binom n2}},
\]
as asserted.
\end{proof}

\begin{corollary}\label{cor:paired-period}
With the notation of Lemma~\ref{lem:HMNcompact},
\[
 \calP_\lambda^\varepsilon
       \bigl(f_{B_n,\lambda,0}(S_{\lambda})\bigr)
 =\frac{\pi^{n/2}}{\Gamma(n/2)}.
\]
\end{corollary}
\begin{proof}
Apply Corollary~\ref{cor:twist} and
\eqref{eq:HMN-average}; the two factors
$(-1)^{\varepsilon {\binom n2}}$ cancel.
\end{proof}

\begin{theorem}
\label{thm:main}
Let $n\ge2$ and let $\mu$ be a dominant integral $n$-tuple.
The period $J_n(\mu)$ defined in
\eqref{eq:intro-J} satisfies
\begin{equation}\label{eq:main-FR}
  J_n(\mu)
    =\frac{\pi^{n/2}}{\Gamma(n/2)}.
\end{equation}
\end{theorem}
\begin{proof}
Take $\lambda=\lambda_{n,\mu}$ and
$\varepsilon\equiv n-1\pmod2$.
Then $J_n(\mu)
=\calP_\lambda^\varepsilon\bigl(f_{B_n,\lambda,0}(S_{\lambda})\bigr)$.
The assertion follows from Corollary~\ref{cor:paired-period}.
\end{proof}

\begin{remark}\label{rem:HMN-character}
To express the period in the Whittaker convention
used by Hara--Miyazaki--Namikawa, put
$ d_n=\diag(i^{n-1},i^{n-2},\ldots,i,1)$ and define
\[
 J_n^+(\mu)=
 \int_{N_{n-1}(\R)\backslash\GL_{n-1}(\R)}
 \calJ_{\psi_+}\bigl(f_{B_n,\lambda,0}(S_\lambda)\bigr)
       \bigl(d_n\diag(h,1)\bigr)
 \sgn(\det h)^{n-1}\dd\bar h.
\]
For any $f\in I_\lambda$ and $g\in G_n$, the substitution
$u=d_nu'd_n^{-1}$ in the Jacquet integral gives
\[
 \calJ_{\psi_+}f(d_ng)
 =i^{\sum_{j=1}^n(j-1)\lambda_j}
  \calJ_{\psi}f(g).
\]
Indeed, $\psi_+(d_nu'd_n^{-1})=\psi(u')$
and $w_nd_nw_n^{-1}=\diag(1,i,\ldots,i^{n-1})$.
Consequently,
\[
 J_n^+(\mu)
 =i^{\sum_{j=1}^n(j-1)\lambda_{n,\mu,j}}
   J_n(\mu)
 =i^{\sum_{j=1}^n(j-1)\lambda_{n,\mu,j}}
   \frac{\pi^{n/2}}{\Gamma(n/2)}.
\]
  While the normalization of Whittaker functions used by Hara--Miyazaki--Namikawa is
   \[
W_{\lambda,0}^{(+)}(v)
=
(-i)^{\sum_{j=1}^n(j-1)\lambda_j}
\Gamma_n(0;\lambda)\,
\mathcal{J}_{\psi_+}\bigl(f_{B_n,\lambda,0}(v)\bigr),
\]
where $
\Gamma_n(0;\lambda)
=
\prod_{1\leq i<j\leq n}
\Gamma\left(1+\frac{\lambda_i-\lambda_j}{2}\right)
$; see \cite[(2.55)]{HMN} for this notation. We thus have 
\[
\int_{N_{n-1}(\mathbb R)\backslash
\mathrm{GL}_{n-1}(\mathbb R)}
W_{d_\mu,0}^{(+)}(S_{\lambda})
\left(
d_n\operatorname{diag}(h,1)
\right)\cdot
\operatorname{sgn}(\det h)^{n-1}\,d\bar h
=
\Gamma_n(0;d_\mu)
\frac{\pi^{n/2}}{\Gamma(n/2)}.
\]
\end{remark}

\begin{corollary}
\label{cor:essentially-tempered}
Let $\mu=(\mu_1\ge\cdots\ge\mu_n)\in\mathbb Z^n$
and $w\in\mathbb Z$. Put
\[
\lambda_0=2\mu^\vee+2\rho_n,\qquad
d=\lambda_0-w(1,\ldots,1),\qquad
\nu_w=(w/2,\ldots,w/2),
\]
and let $\varepsilon\in\{0,1\}$ satisfy
$\varepsilon\equiv n-1-w\pmod 2$.
For $f\in I_{d,\nu_w}=
 \Ind_{B_n(\C)}^{G_n}
 \left(\bigotimes_{j=1}^n
     \chi^{d_j}\abs{\cdot}_{\C}^{\,w/2}\right)$, define
\[
\begin{aligned}
\mathcal P_{d,w}^{\varepsilon}(f)
&=
\int_{N_{n-1}(\mathbb R)\backslash
\mathrm{GL}_{n-1}(\mathbb R)}
\mathcal J_\psi f(\operatorname{diag}(h,1))
|\det h|^{-w}
\operatorname{sgn}(\det h)^\varepsilon
\,d\bar h,\\
\mathcal C_{d,w}^{\varepsilon}(f)
&=
\int_{O(n)}f(k)\det(k)^\varepsilon\,dk.
\end{aligned}
\]
Then $\mathcal P_{d,w}^{\varepsilon}$ converges and
\[
\mathcal P_{d,w}^{\varepsilon}(f)
=
(-1)^{\varepsilon\binom n2}
\frac{\pi^{n/2}}{\Gamma(n/2)}
\mathcal C_{d,w}^{\varepsilon}(f).
\]
Moreover,
\[
\mathcal P_{d,w}^{\varepsilon}
\in
\operatorname{Hom}_{H_n}
\bigl(
I_{d,\nu_w},
|\det|^w\operatorname{sgn}^{\varepsilon}
\circ\det
\bigr).
\]
In particular, let $S_{d}\in V_d$
be the distinguished vector defined in Lemma~\ref{lem:HMNcompact}. Then, using Hara--Miyazaki--Namikawa's normalized Whittaker model,
\[\int_{N_{n-1}(\mathbb R)\backslash
\mathrm{GL}_{n-1}(\mathbb R)}
W_{d,\nu_w}^{(+)}(S_{d})
\bigl(d_n\operatorname{diag}(h,1)\bigr)
\cdot
|\det h|^{-w}
\operatorname{sgn}(\det h)^\varepsilon
\,d\bar h
=
\Gamma_n(\nu_w;d)
\frac{\pi^{n/2}}{\Gamma(n/2)}.
\]
Here $
 \Gamma_n(\nu_w;d)
   =\prod_{1\le i<j\le n}
      \Gamma\left(1+\frac{d_i-d_j}{2}\right), W_{d,\nu_w}^{(+)}
=
(-i)^{q(d)}\Gamma_n(\nu_w;d)
\mathcal J_{\psi_+}\circ f_{B_n,d,\nu_w},
$ 
In particular, $
\Gamma_n(\nu_w;d)
=
\prod_{i<j}\Gamma\left(1+\frac{d_i-d_j}{2}\right)
=
\Gamma_n(0;\lambda_0)
$ is independent of $w$.
\end{corollary}

\begin{proof}
For $f\in I_{d,\nu_w}$, put
$f^\flat(g)=|\det g|^{-w}f(g)$.
Then $f^\flat\in I_d$ and
\[
\mathcal J_\psi f(g)=|\det g|^w\mathcal J_\psi f^\flat(g).
\]
Consequently,
$
\mathcal P_{d,w}^{\varepsilon}(f)
=\mathcal P_d^\varepsilon(f^\flat),
\mathcal C_{d,w}^{\varepsilon}(f)
=\mathcal C_d^\varepsilon(f^\flat).
$
The first identity to prove then follows from Corollary~\ref{cor:twist},
and the equivariance follows by twisting
the corresponding equivariance of
$\mathcal P_d^\varepsilon$.

By \cite[(2.48)]{HMN}, the restriction of
$f_{B_n,d,\nu_w}(v)$ to $U(n)$ is
$(\tau_d(k)v,E_d)_d$. Lemma \ref{lem:HMNcompact} gives
$
\mathcal C_{d,w}^{\varepsilon}
\bigl(f_{B_n,d,\nu_w}(S_{d})\bigr)
=
(-1)^{\varepsilon\binom n2}.
$
Thus
\[
\mathcal P_{d,w}^{\varepsilon}
\bigl(f_{B_n,d,\nu_w}(S_d)\bigr)
=
\frac{\pi^{n/2}}{\Gamma(n/2)}.
\]
Finally, the change-of-character identity
in Remark~\ref{rem:HMN-character} and HMN's normalization show that the two phase factors cancel.
This proves the asserted formula.
\end{proof}

\end{document}